\documentclass[12pt,reqno]{amsart}
\usepackage{amsfonts}
\usepackage{latexsym,amsmath,amsthm,amssymb,amsfonts}

\numberwithin{equation}{section}
\allowdisplaybreaks

\def\endproof{$\hfill\Box$\\}
\def\s{\,\,\,\,}
\def\R{\mathbb{R}}

\numberwithin{equation}{section}
\newtheorem{theorem}{Theorem}[section]
\newtheorem{lem}[theorem]{Lemma}
\newtheorem{thm}[theorem]{Theorem}

\newtheorem{cor}[theorem]{Corollary}

\def\lan{\langle}
\def\ran{\rangle}

\begin{document}
\title[Radially symmetric Solutions To The Graphic Willmore surface Equation]
{Radially symmetric Solutions To The Graphic Willmore Surface Equation}
\author[J. Chen \& Y. Li]
{Jingyi Chen and Yuxiang Li}
\address{ Department of Mathematics\\ The University of British Columbia, Vancouver, BC V6T1Z2, Canada}
\email{jychen@math.ubc.ca}
\address{Department of Mathematical Sciences, Tsinghua University, Beijing 100084, China}
\email{yxli@math.tsinghua.edu.cn}
\thanks{The first author acknowledges the partial support from NSERC. The second author is partially supported by NSFC. This work was done when the second author visits UBC in the fall of 2014; he is grateful to the Department of Mathematics at UBC and PIMS for providing him a nice research environment. }
\date{October 31, 2014}
\maketitle

\begin{abstract}
We show that a smooth radially symmetric solution $u$ to the graphic Willmore surface equation is either a constant or the defining function of a half sphere in ${\mathbb R}^3$. In particular, radially symmetric entire Willmore graphs in ${\mathbb R}^3$ must be flat. When $u$ is a smooth radial solution over a punctured disk $D(\rho)\backslash\{0\}$ and is in $C^1(D(\rho))$, we show that there exist a constant $\lambda$ and a function $\beta$ in $C^0(D(\rho))$ such that  $u''(r) =\frac{\lambda}{2}\log r+\beta(r)$; moreover, the graph of $u$ is contained in a graphical region of an inverted catenoid which is uniquely determined by $\lambda$ and $\beta(0)$. It is also shown that a radial solution on the punctured disk extends to a $C^1$ function on the disk when the mean curvature is square integrable. 

\end{abstract}



\section{Introduction}

A Willmore surface in ${\mathbb R}^3$ is a smoothly immersed surface that satisfies the equation
\begin{equation}\label{Willmore}
\Delta_g H + \frac{1}{2} H^3-2H^2K=0
\end{equation}
where $H = \kappa_1+\kappa_2$ and $K=\kappa_1\kappa_2$ are the mean curvature and the Gauss curvature of the surface respectively and $\Delta_g$ is the Laplace-Beltrami operator in the metric $g$ on the surface, which is induced from the Euclidean metric in ${\mathbb R}^3$ via the immersion.

For a smooth graph $(x,y,u(x,y))$ in ${\mathbb R}^3$, the area element is $v=\sqrt{1+|Du|^2}$ and $H= \mbox{div}\left(\frac{D u}{v}\right)$, where $D$ is the Euclidean gradient operator on functions in the $xy$-plane. Fro graphs, as shown in \cite{D-D},  the Willmore surface equation \eqref{Willmore} can be written in an elegant form:
\begin{equation}\label{Willmore equation}
\mbox{div}\left(\frac{1}{v}\left(\left(I-\frac{D u\otimes D u}{v^2}\right)D(vH)-\frac{1}{2}H^2D u\right)\right)=0.
\end{equation}

Inspired by the classical theorem of Bernstein for the minimal surface equation, one asks when an entire Willmore graph is a plane. An entire smooth solution to \eqref{Willmore equation} is shown to be affine if $H$ is square integrable \cite{C-L}, \cite{L-S}; the approach taken in \cite{C-L} is geometric in nature, rather than dealing with the fourth order partial differential equation directly. It has remained an interesting open question on the radially symmetric case: is an entire radial solution to \eqref{Willmore equation} constant? In this paper, we answer this question affirmatively.

 \begin{thm}\label{main theorem}
Let $u(r)$ be a smooth solution to \eqref{Willmore equation} defined on a disk $D(\rho)$ in the $xy$-plane centered at the origin with radius
$\rho$, where $r=\sqrt{x^2+y^2}$. Suppose that $u(0)=c$ and the mean curvature of of the graph $(x,y,u(r))$ at $(0,0,c)$ is $2a$. Then
\begin{enumerate}
\item  if $a=0$, then $u\equiv c$;
\item if $a\neq 0$,  the graph is contained in a half sphere centered at $(0,0,c)$ of radius $\frac{1}{|a|}$.
\end{enumerate}
\end{thm}
Consequently, a global radially symmetric solution to \eqref{Willmore equation} must be constant.

The main observation in this paper can be summarized as follows. The fourth order equation \eqref{Willmore equation} reduces to a third order equation due to its divergence free structure; the latter further reduces, in the radially symmetric case, to an ordinary differential equation of second order for $u'$ since the zeroth order term $u$ itself is absent from the equation \eqref{Willmore equation}. The singularity at $r=0$ in the second order equation is surprisingly mild. Even more, the difference of any two solutions to the second order equation satisfies an equation which admits a contraction property that leads to uniqueness of solutions to an initial value problem. The desired rigidity then follows. In section \ref{smooth} we will prove Theorem \ref{main theorem}.

\vspace{.1cm}

In section \ref{singular}, we will extend the techniques developed in section \ref{smooth} to study radially symmetric solutions with a point singularity at the origin. Our main focus is to classify the radially symmetric solutions to the Willmore surface
equation on a punctured disk $D(\rho)\backslash\{0\}$, which is in $C^1(D(\rho))$. Such a solution $u$ will be shown to have the following decomposition property:
$$u''(r)=\frac{\lambda}{2}\log r+\beta(r),\s \mbox{where }\beta\in C^0[0,\rho);$$
which is equivalent to the mean curvature decomposes as:
$$H(r)=\frac{\lambda}{2}\log r+\gamma(r),\s \mbox{where }\gamma\in C^0[0,\rho).$$
Further, we will prove that such a surface is contained in the image of an inversion in ${\mathbb R}^3$ of a catenoid, where the inversion and the catenoid are uniquely determined by the information at the origin: $\lambda$ and $\beta(0)$ (or equivalently $\gamma(0)$). We now describe the corresponding catenoid and  inversion.

A catenoid $\mathcal{C}_c$ may be defined by the parametrization
$$
{\mathcal F}_c(t,\theta)=(|c| \cosh t\cos\theta,|c|\cosh t\sin\theta,ct),
$$
where $c\neq 0$ (here we allow $c<0$).  The inversion of this catenoid under the mapping
$$
I_a(x,y,z)=\frac{(x,y,z-a)}{|(x,y,z-a)|^2}
$$
 is a closed surface with $0\in{\mathbb R}^3$ as its only nonsmooth point, in particular, away from finitely many circles that are perpendicular to the $z$-axis, the inverted catenoid can be written as several graphs over the $xy$-plane, among them, we will identify which one contains the graph of $u$ (translated by the constant vector $(0,0,u(0))$).
Let
$$
h(t)=\frac{|c|\cosh t}{c^2\cosh^2t+(ct-a)^2}
$$
be the distance of $I_a({\mathcal F}_c(t,\theta))$ to the $z$-axis. Direct computations shows that $h$ has only finitely many critical points and denote them by $t_1<\cdots<t_m$. Then the inversion of $
\mathcal{C}_c\backslash (\bigcup_i{\mathcal F}_c(\{t_i\}\times S^1))$
consists of the following $m+1$ graphs
$$
I_a\left({\mathcal F}_c((-\infty,t_1)\times S^1)\right),\s I_a\left({\mathcal F}_c((t_1,t_2)\times S^1\right),\cdots,
I_a\left({\mathcal F}_c((t_m,+\infty)\times S^1\right)
$$
where the first one and the last one are graphs over a punctured disk and the rest are graphs over an annulus. Let
$$
t_{c,a}=t_m,\s R_{c,a}=h(t_{c,a}),\s \Sigma_{c,a}=I_a\left({\mathcal F}_c((t_m,+\infty)\times S^1)\right).
$$
The disk $D(R_{c,a})$ in the $xy$-plane is the largest domain
over which $\Sigma_{c,a}$ can be defined (see Section \ref{inverted catenoid} for more details).

\begin{thm}\label{main2}
Let $u(r)$ be a smooth solution to \eqref{Willmore equation} defined on a disk $D(\rho)\backslash\{0\}$ in the $xy$-plane centered at the origin with radius $\rho$. Suppose  $u\in C^{1}(D(\rho))$. Then
\begin{enumerate}
\item[(1)] there exist a constant $\lambda$ and a function $\beta\in C^0([0,\rho))$ such that
$$
u''(r) = \frac{\lambda}{2}\log r +\beta(r),\,\,\,\,r>0;
$$
\item[(2)] the graph of $u$ is contained in $(0,0,u(0))+\Sigma_{c,a}$ where
$$
c=-\frac{\lambda}{4} \,\,\,\, {and} \,\,\,\,
a= \frac{\lambda}{4}\log\frac{|\lambda|}{8}+\frac{3\lambda}{8}-\frac{\beta(0)}{2}
$$
and $\rho$ is no larger than $R_{c,a}$ and $u\in C^{1,\alpha}(D(\rho))$ for any $0<\alpha<1$.
\end{enumerate}
\end{thm}

The regularity assumption $u\in C^1(D(\rho))$ in Theorem \ref{main2} will be  shown  to hold when the mean curvature $H$ is square integrable over $D(\rho)$.
\begin{cor}\label{main3}
Let $u(r)$ be a smooth solution to \eqref{Willmore equation} defined on a punctured disk $D(\rho)\backslash\{0\}$ in the $xy$-plane. Suppose that the Willmore energy $\int_{
D(\rho)}|H|^2d\mu$ is finite. Then $u$ can be extended at $0$ as a $C^{1,\alpha}$ function on $D(\rho)$ and (1) and (2) in Theorem \ref{main2} hold.
\end{cor}

For rigidity of Willmore surfaces in $\R^3$ that are not necessarily graphs, Bryant showed \cite{Bryant}, among other things, that an embedded closed Willmore surface of genus zero must be a round sphere; Rigoli \cite{R}  showed that for a complete
Willmore immersion $f:\Sigma\rightarrow \R^3$,  if there
exists a vector $a\in \R^3$ such that $\lan n+Hf, a\ran\equiv 0$, then $f(\Sigma)$ must be a plane.

\section{Smooth Radially symmetric solutions}\label{smooth}
Let $f$ be a smooth radially symmetric function defined on the disk
$D(\rho)$ in ${\mathbb R}^2$ centered at the origin with radius $\rho\in (0,\infty]$. The divergence of the smooth radial vector field $X=f(r)\frac{\partial}{\partial r}$, where $r=\sqrt{x^2+y^2}$ for $(x,y)\in{\mathbb R}^2$, is given by
\begin{equation}\label{radial vector}
\mbox{div} X=\frac{(rf)'}{r}=f'+\frac{f}{r}.
\end{equation}
Moreover, $f'(0)$ must be 0. In fact, if we set $\tilde f(x,y)=f(r)$, then $\tilde f(x,y)=\tilde f(-x,-y)$. It follows that both $\tilde f_x$ and $\tilde f_y$ vanish at the origin; hence $$f'(0)=\tilde f_x(0,0)\cos\theta+\tilde f_y(0,0)\sin\theta=0. $$


\vspace{0.2cm}

Now consider the graph defined by a smooth radially symmetric function $u(r)$ on $D(\rho)$. Letting $w=u'$, we have
$$
v=\sqrt{1+w^2}\ \ \ \mbox{and} \ \ \ v'=\frac{ww'}{v}.
$$
Note that
$$
w(0)=u'(0)=0.
$$
We calculate the mean curvature of the graph as follows
\begin{eqnarray}\label{H}
H&=& \mbox{div} \left(\frac{w}{v}\frac{\partial}{\partial r} \right) \\
&=&\left( \frac{w}{v} \right)'+\frac{w}{rv} \nonumber \\
&=&\frac{w'}{v}-\frac{w^2w'}{v^3}+\frac{w}{rv} \nonumber \\
&=&\frac{w'}{v^3}+\frac{w}{rv}. \nonumber
\end{eqnarray}
Two consequences are immediate. First, at the origin, we have
\begin{equation}\label{H(0)}
H(0)=w'(0)+\lim_{r\rightarrow 0^+}\frac{w}{rv}=2w'(0).
\end{equation}
Second, we have
\begin{equation}\label{(vH)'}
(vH)'=\frac{w''}{v^2}-\frac{2w(w')^2}{v^4}+\frac{w'}{r}-\frac{w}{r^2}.
\end{equation}

To simplify the divergence free radial vector field in the Willmore surface equation \eqref{Willmore equation}, we first compute

\begin{eqnarray}\label{vH}
\left(I-\frac{D u\otimes D u}{v^2}\right)D(vH)&=&\left(\left(1-\frac{u_x^2}{v^2}\right)(vH)_x-\frac{u_xu_y}{v^2}\left(vH\right)_y\right)
\frac{\partial}{\partial x} \\
&&+\left(-\frac{u_xu_y}{v^2}(vH)_x+\left(1-\frac{u_y^2}{v^2}\right)(vH)_y\right)\frac{\partial}{\partial y}\nonumber \\
&=&\left(\frac{1+u_y^2}{v^2}(vH)_x-\frac{u_xu_y}{v^2}(vH)_y\right)
\frac{\partial}{\partial x}  \nonumber \\
&&+\left(-\frac{u_xu_y}{v^2}(vH)_x+\frac{1+u_x^2}{v^2}(vH)_y\right)\frac{\partial}{\partial y}. \nonumber
\end{eqnarray}
Next, we convert the partial derivatives in $x,y$ of the radially symmetric  functions $u,vH$ into the derivatives in $r$:
\begin{equation}\label{r-derivatives}
u_x=u_r\frac{x}{r},\s u_y=u_r\frac{y}{r},\s  (vH)_x=(vH)_r\frac{x}{r},\s (vH)_y=(vH)_r\frac{y}{r}.
\end{equation}
It follows by substituting \eqref{r-derivatives} into \eqref{vH} that

\begin{eqnarray}\label{vH in r}
\hspace{.75cm}\left(I-\frac{D u\otimes D u}{v^2}\right)D(vH)&=&\frac{1}{v^2}\left(\left({1+u_r^2\frac{y^2}{r^2}}\right)(vH)_r\frac{x}{r}-{u_r^2}(vH)_r\frac{xy^2}{r^3}\right)
\frac{\partial}{\partial x} \\
&&+\frac{1}{v^2}\left(-{u_r^2}(vH)_r\frac{x^2y}{r^3}
+\left(1+u_r^2\frac{x^2}{r^2}\right)(vH)_r\frac{y}{r}\right)\frac{\partial}{\partial y}
 \nonumber \\
&=&\frac{1}{v^2}(vH)_r\left(\frac{x}{r}\frac{\partial}{\partial x}+\frac{y}{r}\frac{\partial}{\partial y}\right) \nonumber \\
&=&\frac{1}{v^2}(vH)_r\frac{\partial}{\partial r}. \nonumber
\end{eqnarray}

Define a radially symmetric function by
\begin{equation}\label{def f}
f=\left<\frac{1}{v}\left(\left(I-\frac{D u\otimes D u}{v^2}\right)D(vH)-\frac{1}{2}H^2D u\right),\frac{\partial}{\partial r}\right>.
\end{equation}
Using \eqref{H}, \eqref{(vH)'} and \eqref{vH in r}, we find
\begin{eqnarray}\label{f}
f&=&\frac{1}{v}\left(\frac{1}{v^2}(vH)_r-\frac{1}{2}H^2 w\right)\\
&=&\frac{1}{v}\left(\frac{1}{v^2}\left(\frac{w''}{v^2}-\frac{2w(w')^2}{v^4}+\frac{w'}{r}-\frac{w}{r^2}\right)\right)-\frac{1}{2v}H^2 w \nonumber\\
&=&\frac{w''}{v^5}-\frac{2w(w')^2}{v^7}+\frac{1}{v^3}\left(\frac{w'}{r}-\frac{w}{r^2}\right)
-\frac{w}{2v}\left(\frac{(w')^2}{v^6}+\frac{w^2(w')^2}{r^2v^2}+2\frac{w'w}{rv^4}\right) \nonumber\\
&=&\frac{1}{v^5}\left(w''+\left(-2\frac{w}{v^2}-\frac{w}{2v^2}\right)(w')^2
+\left(v^2-w^2\right)\frac{w'}{r}-v^2\frac{w}{r^2}-v^2\frac{w^3}{2r^2}\right) \nonumber\\
&=&\frac{1}{v^5}\left(w''-\frac{5w}{2\left(1+w^2\right)}(w')^2+\frac{w'}{r}-\left(1+w^2\right)\frac{w}{r^2}-\left(1+w^2\right)\frac{w^3}{2r^2}\right) \nonumber\\
&=&\frac{1}{v^5}\left(w''+\left(\frac{w}{r}\right)'-\frac{5w}{2\left(1+w^2\right)}(w')^2-\frac{w^3}{2r^2}\left(3+w^2\right)\right).\nonumber
\end{eqnarray}

\begin{lem}\label{lemma 1}
Let $u$ be a smooth radially symmetric function that solves the Willmore surface equation \eqref{Willmore equation} on a disk $D(\rho)$ in the $xy$-plane centred at the origin.  Then $f=0$.
\end{lem}

\proof Using the function defined in \eqref{def f}, the Willmore surface equation \eqref{Willmore equation} reads
\begin{equation}\label{radial W equation}
\mbox{div} \left(  f \frac{\partial}{\partial r}  \right) = 0.
\end{equation}
By \eqref{radial vector}, the above equation reduces to
$$
(rf)'=0.
$$
Then
\begin{equation}\label{f equation}
f=\frac{C}{r}
\end{equation}
for some constant $C$. Furthermore, for any $0<r_0<\rho$, integrating the Willmore equation \eqref{radial W equation} and using Stokes' theorem
$$
0=\int_{D(r_0)}\mbox{div}\left( f\frac{\partial }{\partial r}\right)=\int_{\partial D(r_0)}f=2\pi  \,C.$$
Thus, $C=0$. \endproof

Then, Lemma \ref{lemma 1} and \eqref{f} yield

\begin{lem}\label{lemma 2}
$F=(x,y,u(r))$ is Willmore if and only if
\begin{equation}\label{equation.w}
w''+\left(\frac{w}{r}\right)'-\varphi(\omega)=0,
\end{equation}
where
$$\varphi(\omega)=\frac{5w}{2(1+w^2)}(w')^2+\frac{w^3}{2r^2}(3+w^2).$$
\end{lem}

For \eqref{equation.w}, we now establish a uniqueness result for an initial value problem.
\begin{lem}\label{uniqueness}
There exists $\epsilon>0$ such that the equation \eqref{equation.w} has at most one solution on $[0,\epsilon]$ with
initial data $w(0)=0$, and $w'(0)=a$.
\end{lem}

\proof Let $w_1$ and $w_2$ be solutions to \eqref{equation.w} on $[0,r_0]$ for some $r_0>0$. Then we have
\begin{equation}\label{equation.w1-w2}
(w_1-w_2)''+(\frac{w_1-w_2}{r})'=\varphi(w_1)-\varphi(w_2).
\end{equation}
We have
\begin{eqnarray}\label{varphi(w1)-varphi(w2)}
\varphi(w_1)-\varphi(w_2)&=&\frac{5w_1}{2\left(1+w_1^2\right)}(w_1')^2+\frac{w_1^3}{2r^2}\left(3+w_1^2\right)
-\frac{5w_2}{2\left(1+w_2^2\right)}(w_2')^2-\frac{w_2^3}{2r^2}\left(3+w_2^2\right) \nonumber\\
&=&\left(\frac{5w_1}{2(1+w_1^2)}-\frac{5w_2}{2\left(1+w_2^2\right)}\right)(w_1')^2
+\frac{5w_2(w_1'+w_2')}{2(1+w_2^2)}\left(w_1'-w_2'\right)\nonumber \\
&&+3\,\frac{w_1^3-w_2^3}{2r^2}+\frac{w_1^5-w_2^5}{2r^2}.
\end{eqnarray}

To estimate the three terms on the right hand side of the above identity, we first observe
$$
\left|\frac{w_1}{(1+w_1^2)}-\frac{w_2}{(1+w_2^2)}\right| =
\left|\frac{(1-w_1w_2)(w_1-w_2)}{(1+w_1^2)(1+w_2^2)}\right|\leq |w_1-w_2|.
$$
Next, let
\begin{equation}\label{lambda}
\Lambda=\|w_1\|_{C^1([0,r_0])}+\|w_2\|_{C^1([0,r_0])}.
\end{equation}
Since $w_i(0)=0$, for $r\leq r_0$, it holds
$$
|w_i(r)| = \left| \int^r_0 w'(s)ds\right| \leq r \|w_i\|_{C^1([0,r_0])}.
$$
There exists some positive constant $C(\Lambda)$ which depends only on $\Lambda$ such that:
\begin{eqnarray*}
\left|\frac{w_i}{r}\right|&\leq& C(\Lambda),\s i=1,2;\\
\left|\frac{w_1^3-w_2^3}{r^2}\right| &=& \frac{|w_1^2+w_1w_2+w_2^2|}{r^2}|w_1-w_2| \leq C(\Lambda)|w_1-w_2|, \\
\left|\frac{w_1^5-w_2^5}{r^2}\right| &\leq& C(\Lambda)r^2|w_1-w_2|.
\end{eqnarray*}
Then
\begin{equation}\label{diff}
|\varphi(w_1)-\varphi(w_2)|\leq C(r|w_1'-w_2'|+|w_1-w_2|),\s r\in [0,r_0],
\end{equation}
where $C$ depends only on $\Lambda$ and $r_0$.

Let $\phi=w_1-w_2$. We have $\phi(0)=\phi'(0)=0$, hence $\lim_{r\rightarrow 0^+}
\frac{\phi}{r}=0$; therefore by \eqref{equation.w1-w2},
$$
\phi'(r)+\frac{\phi(r)}{r}=\int_0^r(\varphi(w_1(t))-\varphi(w_2(t)))dt.
$$
Then
$$
(r\phi(r))'=r\int_0^r(\varphi(w_1(t))-\varphi(w_2(t)))dt.
$$
By \eqref{diff},
\begin{eqnarray}\label{contraction}
|(r\phi(r))'|&\leq& Cr\int_0^r(|\phi(t)|+|t\phi'(t)|)dt \\
&=&Cr\int_0^r(|\phi(t)|+\left|(t\phi(t))'-\phi(t)\right|)dt \nonumber\\
&\leq&2Cr\int_0^r(|\phi(t)|+|(t\phi(t))'|)dt \nonumber\\
&=& 2Cr\int_0^r\left(\left|\frac{1}{t}\int_0^t(s\phi(s))'ds\right|+|(t\phi(t))'|\right)dt \nonumber\\
&\leq& 4Cr^2\|(t\phi)'\|_{C^0([0,r])}, \nonumber
\end{eqnarray}
where $r\in(0,r_0]$.

Given $r_1\in(0,r_0)\cap (0,\frac{1}{2\sqrt{2C}})$, if $(r\phi)'$ is not identically $0$ on $[0,r_1]$, then we can find $r_2\in (0,r_1]$, such that
$$
0<|(r\phi)'(r_2)|=\|(r\phi)'\|_{C^0([0,r_2])}.
$$
Thus, \eqref{contraction} yields
$$
0<|(r\phi)'(r_2)|\leq 4Cr_2^2|(r\phi)'(r_2)|<\frac{1}{2}|(r\phi)'(r_2)|,
$$
which is impossible. Therefore, $(r\phi)'=0$ on $[0,r_1]$. It follows that $r\phi$ is a constant which must be 0 since $\phi(0)=0$. We conclude that $\phi=0$ on $[0,r_1]$. Then $w_1=w_2$ on $[0,\epsilon]$ for $\epsilon =\min\{r_0, \frac{1}{2\sqrt{2}C}\}$.
\endproof

\noindent {\it Proof of Theorem \ref{main theorem}}: We divide the proof into two cases.

1) When $a=0$,  it is clear that $0$ is solution to \eqref{equation.w}. Since $w=u'$ is also a solution to \eqref{equation.w}, by Lemma \ref{uniqueness}, $w=0$ on $[0,\epsilon]$ for some $\epsilon>0$. Over the interval $[\epsilon,+\infty)$, \eqref{equation.w}
 is a second order ordinary differential equation with smooth coefficients, by the standard uniqueness theorem in ODE,  $w=0$ on $[\epsilon,+\infty)$ hence on $[0,+\infty)$. Then $u$ is a constant.

2)  Let  $a\neq 0$. It is well-known that the half sphere $F=(x,y,-\mbox{sign}(a)\sqrt{\frac{1}{a^2}-r^2})$ is a
 Willmore surface, i.e. $\left(-\mbox{sign}(a)\sqrt{\frac{1}{a^2}-r^2}\right)'$ satisfies the equation \eqref{equation.w}.  Moreover,
$$\left(-\mbox{sign}(a)\sqrt{\frac{1}{a^2}-r^2}\right)''_{r=0}=a.$$
Similar to  case 1), from Lemma \ref{uniqueness}
$$
w(r)= -\left(\mbox{sign}(a)\sqrt{\frac{1}{a^2}-r^2}\right)'
$$
on $[0,\frac{1}{|a|})$, which is the maximum interval on which $w$ is defined.  Then
$$
u=-\mbox{sign}(a)\sqrt{\frac{1}{a^2}-r^2}+c
$$
where $c$ is a constant and $r\in[0,\frac{1}{|a|})$. \endproof

For an entire radial solution to \eqref{Willmore equation}, the second case in Theorem \ref{main theorem} cannot occur. Therefore, we have
\begin{cor}
A smooth radially symmetric Willmore graph $F=(x,y,u(r))$
must be a plane that is parallel to the $xy$-plane.
\end{cor}

\section{Radially symmetric solutions with a point singularity}\label{singular}

In this section, we study radially symmetric Willmore surface which has a
singularity at $(0,0,0)$.

\subsection{Regularity and uniqueness of solutions}
Examining the proof of Lemma \ref{lemma 2}, due to the singularity at 0, we can only conclude from \eqref{f equation} that
$f=\frac{\lambda}{r}$ for some $\lambda\in\R$, while $\lambda =0$ when $u$ is smooth at $0$. Therefore, the equation for $w = u'$ takes a different form when singularity presents. We have

\begin{lem}\label{lemma}\label{singular Willmore equation}
Let $u\in C^\infty(D(\rho)\backslash\{0\})$ be a  radially symmetric function that solves the Willmore surface equation \eqref{Willmore equation} on $D(\rho)\backslash\{0\}$ in the $xy$-plane centred at the origin.  Then $f(r)=\frac{\lambda}{r}$ for some $\lambda\in\R$ where $f$ is defined in \eqref{f}, i.e. $F=(x,y,u(r))$ is Willmore if and only if
\begin{equation}\label{equation.w.with singularity}
w''+\left(\frac{w}{r}\right)'=\varphi(\omega)+\frac{\lambda}{r}v^5,
\end{equation}
where
$$\varphi(\omega)=\frac{5w}{2(1+w^2)}(w')^2+\frac{w^3}{2r^2}(3+w^2).$$
\end{lem}

We now establish a uniqueness result for an initial value problem of \eqref{equation.w.with singularity}.
\begin{lem}\label{uniqueness2}
There exists $\epsilon>0$ such that the equation \eqref{equation.w.with singularity} has at most one solution on $(0,\epsilon]$ with the
initial data
\begin{equation}\label{initial}
\lim_{r\rightarrow 0^+}w(r)=0,\s \lim_{r\rightarrow 0^+}\left(w'(r)-\frac{\lambda}{2}\log r\right)=b.
\end{equation}
\end{lem}

\proof Let $w_1$ and $w_2$ be solutions to \eqref{equation.w.with singularity} on $(0,r_0]$ for some $0<r_0<1$ with \eqref{initial} satisfied. Let $\phi=w_1-w_2$. It follows that
$$
\lim_{r\rightarrow 0^+}\phi'(r)=\lim_{r\rightarrow 0^+}\phi(r)=0
$$
and $\phi$ can be extended to a function in $C^1([0,r_0])$ with
$\phi(0)=0$ and $\phi'(0)=0$.
Hence $\lim_{r\rightarrow 0^+}
\frac{\phi}{r}=0$, and by integration
$$
\phi'(r)+\frac{\phi(r)}{r}=\int_0^r\left(
\varphi(w_1(t))-\varphi(w_2(t))+\lambda\frac{v_1^5(t)-v_2^5(t)}{t}\right)dt.
$$
Then
$$
(r\phi(r))'=r\int_0^r\left(\varphi(w_1(t))-\varphi(w_2(t))+\lambda\frac{(1+w_1^2)^\frac{5}{2}-(1+w_2^2)^\frac{5}{2}}{t}\right)dt.
$$
From \eqref{initial}, we  may find $\Lambda_1=\Lambda_1(r_0,b,\lambda)$ such that
\begin{equation}\label{lambda1}
|w_1'(r)|+|w'_2(r)|<\Lambda_1 \log \frac{1}{r}.
\end{equation}
and
$$
\frac{|w_i(r)|}{r}\leq\frac{1}{r}\int_0^r|w_i'| \leq\Lambda_1 \log \frac{1}{r},\s i=1,2.$$
Moreover,
$$
\left|\frac{(1+w_1^2(r))^\frac{5}{2}-
(1+w_2^2(r))^\frac{5}{2}}{r}\right|\leq \frac{C}{r}|w_1^2-w_2^2|
\leq C|\phi|\log\frac{1}{r}.
$$
Applying \eqref{varphi(w1)-varphi(w2)}, similar to the proof of \eqref{diff} by using \eqref{lambda1} instead of \eqref{lambda}, we have
$$
\left|\varphi(w_1(r))-\varphi(w_2(r))+\lambda\frac{(1+w_1^2(r))^\frac{5}{2}-
(1+w_2^2(r))^\frac{5}{2}}{r}\right|\leq C(|\phi(r)|+r|\phi'|)\log^2\frac{1}{r}.
$$
Then using the arguments in \eqref{contraction}, we have 
\begin{eqnarray}\label{contraction1}
|(r\phi(r))'|&\leq& Cr\int_0^r(|\phi(t)|+|t\phi'(t)|)\log^2\frac{1}{t}dt \\
&\leq& 4Cr\left(r\log^2r-2r\log r+2r\right)\|(t\phi)'\|_{C^0([0,r])}\, \nonumber
\end{eqnarray}
where $r\in(0,r_0]$.

Given $r_1\in (0,r_0)$, such that 
$$
4Cr_1\left(r_1\log^2r_1-2r_1\log r_1+2r_1\right)\leq\frac{1}{2}.
$$ 
If $(r\phi)'$ is not identically $0$ on $[0,r_1]$, then we can find $r_2\in (0,r_1]$, such that
$$
0<|(r\phi)'(r_2)|=\|(r\phi)'\|_{C^0([0,r_2])}.
$$
Thus, \eqref{contraction1} yields
$$
0<|(r\phi)'(r_2)|\leq\frac{1}{2}|(r\phi)'(r_2)|,
$$
which is impossible. Therefore, $(r\phi)'=0$ on $[0,r_1]$, which implies that  $w_1=w_2$ on $[0,\epsilon]$ for $\epsilon =\min\{r_0, r_1\}$.
\endproof

\begin{lem}\label{regularity}
Let  $w\in C^2((0,a_0])$ be a solution to \eqref{equation.w.with singularity} on $(0,a_0]$. If $\lim_{r\rightarrow 0^+}w(r)=0$,
then $w(r)-\frac{\lambda}{2}r\log r\in C^1([0,a_0])$, and consequently $\lim_{r\rightarrow 0^+}\left(w'(r)-\frac{\lambda}{2}\log r\right)$ exists.
\end{lem}

\proof
We assume $|w|<\epsilon_1<1$ on $(0,a], a\leq a_0, a<1$.
First of all, we prove
\begin{equation}\label{energy}
\int_r^a\left((w')^2+\frac{w^2}{2t^2}\right)<C\left(\frac{w^2}{2r}+|w|\log\frac{1}{r}+\int_r^a\frac{|w|}{t}dt+1\right),\s r\in(0,a),
\end{equation}
where $C$ is a positive constant depends only on $a$,  $w'(a)$, $\lambda$ and the upper bound of $w$ in $[0,a]$. For simplicity, we will use $C$ to denote some uniform constant in the rest of the proof.

Multiplying  both sides of \eqref{equation.w.with singularity} by $w$ and then integrating from $r$ to $a$, since $w\varphi\geq 0$, we have
$$
\int_r^{a}\left(w''w+\left(\frac{w}{t}\right)'w\right)\geq \int_r^a\frac{\lambda}{t}v^5w\geq -C\int_r^a\frac{|w|}{t}dt.
$$
Then integrating by parts implies
$$
-w'(r)w(r)-\frac{w^2(r)}{r}-\int_r^a\frac{(w^2)'}{2t}-\int_r^a(w')^2
\geq -C\left(1+\int_r^a\frac{|w|}{t}dt\right),
$$
then applying integration by parts again and absorbing the terms evaluated at $a$ into $C$,  we have
$$
-w'(r)w(r)-\frac{w^2(r)}{r}+\frac{w^2(r)}{2r}-\int_r^a\frac{w^2}{2t^2}-\int_r^a(w')^2
\geq -C\left(1+\int_r^a\frac{|w|}{t}dt\right).
$$
Rearranging terms, we arrive at
\begin{equation}\label{energy2}
\int_r^a\left((w')^2+\frac{w^2}{2t^2}\right)dt\leq C\left(1+\int_r^a\frac{|w|}{t}dt\right)-w(r)\frac{(rw(r))'}{r}+\frac{w^2(r)}{2r}.
\end{equation}
To estimate $(rw(r))'$, integrating \eqref{equation.w.with singularity} again,
\begin{equation}\label{(rw)'}
\frac{(rw(r))'}{r}=-\int_r^a\left(\frac{5w(w')^2}{2(1+w^2)}+
\frac{3w^3+w^5}{2t^2}+\lambda\frac{v^5-1}{t}\right)dt
+\lambda \log\frac{r}{a}+w'(a)+\frac{w(a)}{a}.
\end{equation}
Then, absorbing the terms evaluated at $a$ into a uniform constant $C$, we have 
\begin{equation}\label{(rw)'2}
\left|\frac{(rw(r))'}{r}\right|\leq
\epsilon_1 \int_r^a\left((w')^2+\frac{w^2}{2t^2}\right)+C\log\frac{1}{r}.
\end{equation}

Putting \eqref{(rw)'2} into \eqref{energy2} and recalling $|w|<\epsilon_1$, we get
$$
\int_r^a\left((w')^2+\frac{w^2}{2t^2}\right)\leq C\left(1+\int_r^a\frac{|w|}{t}dt\right)+\epsilon_1^2\int_r^a\left((w')^2+\frac{w^2}{2t^2}\right)+C|w|\log\frac{1}{r}+\frac{w^2(r)}{2r}
$$
and in turn we conclude the proof of \eqref{energy}.

Secondly, we prove
\begin{equation}\label{w}
|w(r)|<Cr\log \frac{1}{r}.
\end{equation}
By \eqref{energy} and \eqref{(rw)'2},
$$|rw'+w|\leq C\left(w^2+r\log\frac{1}{r}\right).$$
Fix an $\epsilon_2\in (0,1)$ and choose $a'$ such that $|w|<{\epsilon_2}/{C}$ on
$(0,a')$. It is well-known that
$$
\big|\,|rw|'\,\big|\leq |(rw)'|.
$$
Therefore
$$
\left|(r|w|)'\right|\leq \epsilon_2|w|+Cr\log\frac{1}{r},\s r\in(0,a').
$$
Then
$$r|w|'+(1-\epsilon_2)|w|\leq Cr\log\frac{1}{r}.$$
Thus
$$
\left(r^{1-\epsilon_2}|w|\right)'\leq Cr^{1-\epsilon_2}\log\frac{1}{r},
$$
then
$$
r^{1-\epsilon_2}|w(r)|\leq C\int_0^rt^{1-\epsilon_2}\log\frac{1}{t}dt\leq \frac{C}{2-\epsilon_2}
\left(r^{2-\epsilon_2}\log\frac{1}{r}+r^{2-\epsilon_2}\right)
$$
and \eqref{w} is established by dividing $r^{1-\epsilon_2}$.

Now we can finish the proof.
Putting \eqref{w} into \eqref{energy}, we are led to
$$
\int_0^a\left((w')^2+\frac{w^2}{2t^2}\right)dt<+\infty,
$$
then
\begin{equation}\label{finite}
\int_0^a\left|\frac{5w}{2(1+w^2)}(w')^2+\frac{3w^3+w^5}{2t^2}+\lambda\frac{v^5-1}{t}
\right|dt<+\infty.
\end{equation}
For convenience, set
$$
\psi(r)=-\int_r^a\left(\frac{5w}{2(1+w^2)}(w')^2+\frac{3w^3+w^5}{2t^2}+\lambda\frac{v^5-1}{t}\right)dt,
$$
which is continuous in $[0,a]$ by \eqref{finite}. By \eqref{(rw)'},
$$
w(r)=\frac{1}{r}\int_0^rt\psi(t)dt+\frac{\lambda}{2}r\log r+Br,
$$
where 
$$ 
B=\frac{1}{2}\left(-\lambda\log a+w'(a)+\frac{w(a)}{a}\right)-\frac{\lambda}{4}.
$$
Thus,
$$
\lim_{r\to 0^+}\left(w(r)-\frac{\lambda}{2}\log r\right)'=B+\frac{\psi(0)}{2}.
$$
It is easy to check that the function $w(r)-\frac{\lambda}{2}r\log r$ is differentiable at 0 with derivative $B+\frac{\psi(0)}{2}$. Therefore, $w(r)-\frac{\lambda}{2}r\log r\in C^1([0,a_0])$. 
\endproof



\subsection{Inverted catenoids}\label{inverted catenoid}
A catenoid $\mathcal{C}_c$ is defined by
$${\mathcal F}_c(t,\theta)=(|c| \cosh t\cos\theta,|c|\cosh t\sin\theta,ct),$$
where $c\neq 0$.
The inversion of $\mathcal{C}_c$ under the map $I_a:{\mathbb R}^3\to{\mathbb R}^3$ defined in the introduction is given by
$$
(x,y,z)=\left(\frac{|c| \cosh t\cos\theta}{R^2},\frac{|c|\cosh t\sin\theta}{R^2},
\frac{ct-a}{R^2}\right),
$$
where
$$
R^2=c^2\cosh^2t+(ct-a)^2.
$$
Let
$$
h(t)=\frac{|c|\cosh t}{R^2},\s \bar{h} (t)=\frac{ct-a}{R^2}.
$$
Then $I_a(\mathcal{C}_c)$ is the surface of revolution generated by rotating the curve
$$
\sigma(t)=(h(t),0,\bar{h}(t))
$$
about the $z$-axis.
Since $h(+\infty)=h(-\infty)=0$, $h$ has at least 1 critical point. Direct calculation shows that $h$ has only finitely many (in fact, at most three) critical points. A tangent vector of $\sigma$ is given by $\sigma'(t)=(h'(t),0,\bar h'(t))$, so $t$ is a critical point of $h$ if and only if the tangent line of $\sigma(t)$ is parallel to the $z$-axis. Let $t_1<\cdots<t_m$ be the critical points of $h$ and set $t_0=-\infty$, $t_{m+1}=+\infty$, then
$I_a({\mathcal F}_c((t_i,t_{i+1})\times S^1))$ is a graph over
$D(h(t_{i+1}))\backslash D(h(t_i))$, which cannot be extended as a smooth graph over a larger domain, where $i=0, \cdots, m$.

\begin{lem}\label{catenoid}
Given $\lambda,b\in \R$, there exists a unique pair of numbers $c$ and $a$, such that
$\Sigma_{c,a}$ is the graph of a function $U(r)$ with
$$
\lim_{r\rightarrow 0^+}\left(U''(r)-\frac{\lambda}{2}\log r\right)=b.
$$
Moreover, $U'$ satisfies the equation \eqref{equation.w.with singularity}.
\end{lem}

\proof

After changing parametrization, a catenoid can be written as
$$
(\rho\cos\theta,\rho\sin\theta,u(\rho)), \s \mbox{where}\s  u(\rho)=c \operatorname{arccosh}
\frac{\rho}{|c|}.$$
In this parametriation, $R=\sqrt{\rho^2+(u(\rho)-a)^2}$.
Suppose $\Sigma_{c,a}$ is the graph of $U(r)$. We have
$$
r=\frac{\rho}{R^2},\s \mbox{ and }\s U(r)=\frac{u(\rho)-a}{R^2}=\frac{r}{\rho}(u(\rho)-a),\s
r\in (0,R_{c,a}).
$$
We now calculate the expansion of $U$ at $r=0$. It is clear 
\begin{equation}\label{rho}
\frac{1}{\rho}=r\left(1+\frac{(u(\rho)-a)^2}{\rho^2}\right)
\end{equation}
and 
\begin{eqnarray}\label{U}
U(r)&=&r^2(u(\rho)-a)\left(1+\frac{(u(\rho)-a)^2}{\rho^2}\right)\\
     &=&r^2u(\rho)-ar^2+\frac{r^2}{\rho^2}O(u^3(\rho)), \s \mbox{ as }\rho\rightarrow+\infty.\nonumber
\end{eqnarray}
By the definition of $\Sigma_{c,a}$,  $\frac{u(\rho)}{c}\rightarrow+\infty$ as $\rho\rightarrow+\infty$,
then by taking log of 
$$
e^{\frac{u(\rho)}{c}}\frac{1+e^{-2\frac{u(\rho)}{c}}}{2}=\cosh \frac{u(\rho)}{\rho}=\frac{\rho}{|c|},
$$
we are led to 
\begin{eqnarray*}
u(\rho)&=&c\log\rho+c\log 2-c\log|c|-c\log\left(1+e^{-2\frac{u(\rho)}{c}}\right)\\
&=&c\log\rho+c\log 2-c\log|c|+o(1), \s \mbox{as $\rho\to +\infty$.}
\end{eqnarray*}
Using \eqref{rho}, we get
\begin{eqnarray*}
u(\rho)&=&-c\log r+c\log 2-c\log|c|-\log\left(1+\frac{(u(\rho)-a)^2}{\rho^2}\right)+o(1)\\
&=&-c\log r+c\log 2-c\log|c|+o(1),\s\mbox{as $\rho\to +\infty$}.
\end{eqnarray*}
By \eqref{U},
$$
U(r)=-cr^2\log r+(c\log 2-c\log|c|-a)r^2+o(r^2),\s\mbox{ as $r\to 0$.}
$$
It follows that $U(r)+cr^2\log r$ is in $C^2([0,\rho))$ and
$$
\lim_{r\to 0^+}\left(U''(r)+2c\log r\right)= 2(c\log 2-c\log|c|-a)-3c.
$$
Set
$$
c=-\frac{\lambda}{4},\s a=\frac{\lambda}{4}\log\frac{|\lambda|}{8}+\frac{3\lambda}{8}
-\frac{b}{2}.
$$
Then the uniquely determined $\Sigma_{c,a}$ is the graph of $U$.
\endproof

\subsection{Proof of Theorem \ref{main2} and Corollary \ref{main3}} Note that a smooth radially symmetric function $f$ defined on $D(\rho)\backslash\{0\}$ can be extended at 0 as a function in $C^1(D(\rho))$ if and only if $\lim_{r\to 0^+}f'(r)=0$. 

\vspace{.2cm}

\noindent {\it Proof of Theorem \ref{main2}}. Since $u\in C^1(D(\rho))$, $\lim_{r\to 0^+}w(r)=0$.  
By Lemma \ref{regularity}, we can write
$$
w(r)'-\frac{\lambda}{2}\log r = \beta(r)
$$
for some function $\beta\in C^0([0,\rho))$. So
$$
\lim_{r\rightarrow 0}(w'-\frac{\lambda}{2}\log r)=\beta(0).
$$
By Lemma \ref{catenoid}, there exist $c$ and $a$, such that $\Sigma_{c,a}$ is the graph
of a function $U(r)$ with
$$
\lim_{r\rightarrow 0}(U''-\frac{\lambda}{2}\log r)=\beta(0)
$$
and $U'$ satisfies \eqref{equation.w.with singularity}. Then by Lemma \ref{uniqueness2},
$U'=w$ in $[0,\epsilon]$ for some $\epsilon$. Note that the coefficients in \eqref{equation.w.with singularity}
are smooth on $[\epsilon,\rho)$. The standard uniqueness theorem in ODE theory asserts $U'=w$ on $[\epsilon, \rho)$ and consequently $U'=w$ on $[0,\rho)$.
\endproof


\noindent{\it Proof of Corollary \ref{main3}.}  By a direct calculation, the first and the second fundamental forms of $F$ are
$$g=\left(\begin{array}{ll}
                1+w^2&0\\
                0&r^2\end{array}\right),\s
A=\left(\begin{array}{ll}
                \frac{w'}{v}&0\\
                0&\frac{rw}{v}\end{array}\right)$$
respectively.  Then
\begin{equation}\label{A}
|A|^2=\left(\frac{1}{v^2}\frac{w'}{v}\right)^2+\left(\frac{1}{r^2}\frac{rw}{v}\right)^2=
\left(\frac{w'}{v^3}\right)^2+\left(\frac{w}{rv}\right)^2
\end{equation}
and
\begin{eqnarray*}
\int_{D(\rho)\backslash D(\epsilon)}(|A|^2-|H|^2)d\mu
&=&-2\int_{D(\rho)\backslash D(\epsilon)}\frac{w'}{v^3}\frac{w}{rv}d\mu\\
&=&-4\pi\int_\epsilon^{\rho}\frac{ww'}{v^3}\\
&=&4\pi\int_\epsilon^\rho\left(\frac{1}{v}\right)'dr\\
&=&4\pi\left(\frac{1}{v(\rho)}-\frac{1}{v(\epsilon)}\right).
\end{eqnarray*}
Thus,
\begin{equation}\label{A finite}
\int_{D(\rho)}|A|^2<+\infty.
\end{equation}

We now prove $\lim_{r\to 0^+}w(r)=0$. Assume there is a sequence
$r_k\rightarrow 0$, such that $w(r_k)\rightarrow c\neq 0$, where $c$ may be $\infty$. Define
$$
u_k(r)=\frac{u_k(r_kr)-u_k(r_k)}{r_k},\s w_k(r)=u_k'(r),\s F_k=(x,y,u_k),\s \Sigma_k=
F_k\left(D\left(\frac{\rho}{r_k}\right)\backslash \{(0)\}\right).
$$
Obviously, $w_k(1)=w(r_k)\rightarrow c$ as $k\to\infty$.
We may assume that the unit normal vector $n_k(p)$ of $\Sigma_k$ at $p=(1,0,0)$  converges to a unit vector $\alpha$ which is not orthogonal to the $xy$-plane,
i.e.
\begin{equation}\label{not 0}
\alpha\cdot (1,0,0)\neq 0.
\end{equation}
This is because the angle between the tangent plane $T_{p}\Sigma_k$ and the $xy$-plane
is uniformly bounded away from 0 when $c\not=0$.
Thus, we may assume that $\Sigma_k$ is the graph of some function $\tilde{u}_k(y,z)$ near $p$ over the $yz$-plane.  For convenience, we set $\tilde{F}_k=(\tilde{u}_k(y,z),y,z)$ and
$e=\frac{\partial\tilde{F}_k}{\partial y}$. Since $\Sigma_k$ is radially symmetric, the unit circle
$\{(\cos\theta,\sin\theta,0):\theta\in[0,2\pi)\}$ is in each $\Sigma_k$. Then
$\tilde{u}_k(y,0)=\sqrt{1-y^2}$ and
$|e(p)|=1$. For the second fundamental form $A_k$ of $\Sigma_k$ we have
\begin{equation}\label{A_k}
A_k(e(p),e(p))=\frac{\partial^2\tilde{F}_{k}}{\partial y^2}(p)\cdot n_k(p)=
(-1,0,0)\cdot n_k(p).
\end{equation}
However, for any $\delta>0$, we have by \eqref{A finite}
$$
\lim_{k\rightarrow+\infty}\int_{D(\delta)}|A_k|^2d\mu_k=
\lim_{k\rightarrow+\infty}\int_{D(r_k\delta)}|A|^2d\mu=0.
$$
Applying the curvature estimates in Theorem 2.10 in \cite{K-S1},
$$|A_k|(p)\leq C\|A_k\|_{L^2(B^3_1(p))}\leq C\|A_k\|_{L^2(D(2))}\rightarrow 0,$$
where $B^3_1(p)=\{q\in\R^3:|q-p|<1\}$. By \eqref{A_k}, we get
$(1,0,0)\cdot\alpha=0$. This contradicts \eqref{not 0}. Therefore, we have $\lim_{r\to 0^+}w(r)=0$, in turn, $u\in C^1(D(\rho))$. Now the desired results follow from Theorem \ref{main2}. 
\endproof

\end{document}